\documentclass[twocolumn,showpacs,%
  nofootinbib,aps,superscriptaddress,%
  eqsecnum,prd,notitlepage,showkeys,10pt]{revtex4-1}

\usepackage{amssymb}
\usepackage{amsmath}
\usepackage{graphicx}
\usepackage{caption}
\usepackage{subcaption}
\usepackage{dcolumn}
\usepackage{hyperref}
\usepackage{enumerate}
\usepackage{algpseudocode}
\usepackage{url}

\def \N {\mathbb{N}}
\def \Z {\mathbb{Z}}
\def \R {\mathbb{R}}

\def\set#1{\{ #1 \}}
\def \phi {\varphi}

\usepackage{algorithm}
\usepackage{algpseudocode}
\usepackage{amsmath}
\usepackage{amssymb}

\usepackage{setspace}
\usepackage{microtype}
\begin{document}

\title{The Continuous p-Dispersion Problem in Three Dimensions}
\author{Sanjay Manoj}
\affiliation{University of Illinois Urbana-Champaign}
\author{Melkior Ornik}
\affiliation{University of Illinois Urbana-Champaign}

\begin{abstract}
\textbf{Abstract:} The Continuous p-Dispersion Problem (CpDP) with boundary constraints asks for the placement of a fixed number of points in a compact subset of Euclidean space such that the minimum distance between any two points, as well as the points and the boundary of this compact set is maximized. This problem finds applications in facility placement, communication network design, sampling theory, and particle simulation; however, finding optimal solutions is NP-hard and existing algorithms focus on providing approximate solutions in two-dimensional space. In this paper, we introduce an almost-everywhere differentiable optimization model and global optimization algorithm for approximating solutions to the CpDP with boundary constraints in convex and non-convex polyhedra with respect to any metric in a three-dimensional Euclidean space. Our algorithm generalizes two-dimensional dispersion techniques to three dimensions by leveraging orientation, linear-algebraic projections for point-to-face distances, and a ray-casting procedure for point-in-polyhedron testing, enabling optimization in arbitrary convex and non-convex three dimensional polyhedra. We validate the proposed algorithm by comparing with analytical optima where available and empirical benchmarks, observing close agreement with optimal solutions and improvements over empirical benchmarks.
\end{abstract}

\maketitle

\section{Introduction}\label{sec:Intro}

In 1987, \citet{kuby1987} proposed the $p$-Dispersion Problem, asking to locate $p$ facilities on a network so that the minimum separation distance between any pair of open facilities is maximized. Later in 1995, \citet{drezner1995} formalized the $p$-Dispersion Problem, asking for the placement of $p$ points from a subset of $\R^n$ such that the minimum distance between pairs of points is maximized. Subsequently, the Continuous $p$-Dispersion Problem (CpDP) with boundary constraints adds the requirement that a dispersion point must be sufficiently far away from the boundary of the subset of $\R^n$, also often referred to as a container. The $p$-Dispersion Problem is interesting because it formalizes the inherent trade-off between coverage and separation, to find how a limited number of points can be optimally distributed over a domain. The $p$-Dispersion problem is provably NP-complete for the general container \citep{baurfekete2001} in the number of dispersion points $p$. 

Existing work solving the p-Dispersion problem primarily focuses on low-dimensional containers with special geometric properties. For instance, repulsion-based heuristics have been applied to circular containers \citep{graham1998} and non-convex two-dimensional polygons \citep{dai2021}. Other approaches include a 2/3-approximation algorithm for rectilinear polygons \citep{baurfekete2001} and a Tabu Search Global Optimization (TSGO) algorithm for general planar containers \citep{lai2024}. In higher dimensions, work is severely limited: \citet{drezner1995} formulated a nonlinear model restricted to the unit square, while \citet{kazakov2018} developed a grid-based sphere-packing algorithm that struggles with complex, non-convex three-dimensional boundaries.

Compared to previous work, the goal of this paper is to approximate elements of the solution set for the CpDP within general polyhedral containers with respect to any metric in $\mathbb{R}^3$ \cite{lai2024}. The TSGO algorithm proposed by \citet{lai2024} is difficult to adapt to higher dimensions because it is built around two-dimensional logic. Their model relies on measuring distances to the edges of a polygon and uses repeated point-in-polygon tests to check for constraint violations. We extend this approach to solve the CpDP for general polyhedral containers in three dimensions. Our method uses a ray-casting procedure for point-in-container testing and replaces their edge-based distance checks with a formulation that calculates the distance to the faces, edges, and vertices of the three-dimensional boundary. 

Solving the CpDP in three dimensions reflects how spatial dispersion is performed in real-world situations where objects have volume and operate in bounded physical space, such as collision avoidance for autonomous agents \cite{abdulghafoor2021} or transmitter placement in complex indoor environments \cite{loganathan2024}. Importantly, the mathematical abstractions based on \citep{lai2024} established in this work are designed to be extensible to general $n$-dimensional containers. While the primary challenge in such a transition lies in the geometric modeling of arbitrary high-dimensional polyhedra, containers with uniform constraints, such as an $n$-dimensional hypercube with limits restricted between 0 and 1, remain straightforward to model within this framework. We build on the work from \citet{lai2024} by defining a feasibility-residual function to approximately solve the CpDP with boundary constraints for convex and non-convex polyhedral containers in $\R^3$. A feasibility-residual function $E(x) = 0$ if and only if $x$ is feasible solution, and $E(x) > 0$ otherwise. To minimize this function, we utilize the Tabu Search heuristic \cite{glover1997} which improves upon standard local searches by marking previously visited configurations as "tabu" to prevent cycling and encourage the exploration of more optimal regions of the solution space \cite{lai2024}.

We also generalize the Continuous $p$-Dispersion Problem by allowing the distance between points and the boundary to be measured using a metric. This allows the algorithm to handle separation constraints defined by general shapes—like cubes or octahedra—rather than being limited to standard spherical constraints. This change broadens the types of dispersion problems the algorithm can solve and provides a single framework for optimizing point layouts under different geometric requirements. We validate our algorithm's accuracy by comparing its results against optimal results in the unit cube proven by \citet{schaer1966}, measuring the relative deterioration to show how closely our algorithm approaches global maxima. For non-rectangular containers, we benchmark our results against the best-known packings for the unit tetrahedron provided by \citet{kazakov2018}. We also establish new results for non-convex polyhedra and provide visualizations for our best-achieved configurations in several containers.

\textbf{Basic Notation:} We use the following standard notation throughout this paper. Let $d(\cdot, \cdot)$ denote a metric on $\R^n$ following the usual definition \citep{rosen2012}. The distance between an element $x \in \R^3$ and a nonempty subset $X \subseteq \R^3$ is given by $d(x, X) = \inf \set{d(x, y) : y \in X}$. The \textit{open ball} and \textit{closed ball} of radius $r$ centered at $x$ are defined as $B_r^o(x) = \set{y \in \R^3 : d(x, y) < r}$ and $B_r^c(x) = \set{y \in \R^3 : d(x, y) \le r}$, respectively. For any positive integer $n$, let $[n] = \set{1, \dots, n}$. Finally, for a subset $X \subseteq \R^3$, $\partial X$ denotes the boundary of $X$, $X^c = \R^3 \setminus X$ is the complement of $X$, and $X^o = X \setminus \partial X$ is the interior of $X$.

\section{Problem Formulation}

The Continuous $p$-Dispersion Problem (CpDP) asks to place $p$ points, i.e., \textit{sites}, within a compact subset, i.e., a \textit{container}, in $\R^n$ so that the minimum pairwise distance between the sites is maximized \citep{drezner1995}. A variation on the CpDP includes ``boundary constraints'' which asks that the minimum distance between any two points, as well as the points and the boundary of this compact set be maximized \citep{baurfekete2001}. Most existing algorithms that solve the CpDP do so in the setting of bounded, polygonal containers in two dimensions, possibly with holes. In the three dimensions, the natural extension is to consider bounded containers which are polyhedra with polyhedral holes. A \textit{polyhedron} is a closed and bounded solid $C$ whose boundary $\partial C$ can be covered by finitely many planes \cite{berg2008}. A \textit{polyhedral container} is a polyhedron $C$ with holes $U_1, \dots, U_k \subseteq C$ which are also polyhedra.

We denote $C = (V, E, F)$ where $V$ is the set of vertices of $C$, $E$ is the set of edges of $C$ and $F$ is the set of faces of $C$. We require that the interior of a polyhedral container be non-empty, therefore providing the polyhedral container with a consistent orientation with inward normals \citep{munkres1991}. We write $x \in C$ for an element $x \in \R^3$ inside the polyhedron but outside its holes. An element $x \notin C$ is located outside the polyhedron or inside the polyhedral container's holes.

We now state the Continuous $p$-Dispersion Problem (CpDP) \citep{drezner1995} with boundary constraints for polyhedral containers with respect to any metric $d$ on $\R^3$.
\\

\textbf{Problem:} Let $p \in \N$ and $C \subseteq \R^3$ be a polyhedral container. We solve the nonlinear optimization problem:
\begin{align*}
    \max_{c_1, \dots, c_p} \quad & D \\ 
    \text{subject to} \quad & c_1, \dots, c_p \in C, & (i) \\ 
    & d(c_i, c_j) \ge 2D \text{ for all distinct } i, j \in [p], & (ii) \\ 
    & d(c_i, \partial C) \ge D \text{ for all } i \in [p]. & (iii) 
\end{align*}
Note that we will interpret the value of $D$ to be the radius of the ball centered at each dispersion point $c_i$. For a feasible solution $(c_1, \dots, c_p, D) \in (\R^3)^p \times \R$ to the nonlinear optimization problem, we refer to the elements $c_i \in \R^3$ as \textit{dispersion points} and $D \in \R_{> 0}$ as the \textit{radius} of the solution. A \textit{configuration} (of dispersion points) $X$ refers to the set $X = \set{c_1, \dots, c_p}$ for a given feasible solution.

These constraints enforce three geometric requirements: ($i$) the \textit{containment constraint} ensures each point $c_i$ lies within the container; ($ii$) the \textit{dispersion constraint} requires mutual separation, meaning $B_D^o(c_i) \cap B_D^o(c_j) = \emptyset$ for all distinct $i, j \in [p]$; and ($iii$) the \textit{boundary constraint} keeps points away from the edges, ensuring $B_D^o(c_i) \cap \partial C = \emptyset$ for all $i \in [p]$.

From the assumption that the interior of the container $C$ is nontrivial, a feasible solution $(c_1, \dots, c_p, D)$ always exists \citep{dai2021}.

The CpDP is hard to attack naively due to the massive number of non-convex feasible sets of the nonlinear program. As a result, standard nonlinear optimizers often get caught in local minima. Additionally, the number of distance constraints to compute is $\binom{p}2 = \mathcal{O}(p^2)$, and with minimum distance bounds, it is highly likely that a random initial configuration will violate the dispersion constraint (\textit{ii}).

\section{Optimization Model}

In this section we introduce a unified feasibility-residual function used to solve the CpDP with boundary constraints by modeling constraint violation in two parts: \textit{dispersion violation} and \textit{container violation}. Dispersion violation stems from dispersion points being located too close together; namely, when $d(c_i, c_j) < 2D$ for some fixed $D > 0$, constraint (\textit{ii}) is violated. Container violation is a result of dispersion points being located too close to the boundary of the container; namely, when $d(c_i, \partial C) < D$ for some fixed $D > 0$, constraint (\textit{iii}) is violated.

By combining these two components into a single feasibility-residual function, we can address inherent geometric complexities associated with any container. In general, the feasible region for the solution set to the CpDP is non-convex. Moreover, as the number of dispersion points $p$ increases, the feasible region will contain many local optima \citep{lai2024}.

\citet{lai2024} proposed a novel framework for the two-dimensional CpDP by unifying dispersion and boundary constraints into a single, non-negative feasibility-residual function, $E_D$. They conceptualized constraint violations as \textit{energy}, where the total energy of a configuration corresponds to its distance from feasibility. For a configuration of points $X \in (\R^2)^p$ and fixed $D \ge 0$, this function is defined as $\displaystyle{E_D(X) = \sum_{i = 1}^{p - 1} \sum_{j = i + 1}^{p} \ell^2_{ij} + \sum_{i = 1}^p O_i}$ where $\ell_{ij}$ is the dispersion penalty between points $c_i$ and $c_j$, and $O_i$ is the boundary penalty for point $c_i$. A crucial property of this formulation is that $E_D(X) = 0$ if and only if $(X, D)$ is a feasible solution. Because $E_D$ is almost-everywhere differentiable, it allows for the application of efficient gradient-based optimization methods, such as the Limited-memory Broyden-Fletcher-Goldfarb-Shanno (L-BFGS) algorithm \citep{fletcher1987}. This approach proved exceptionally powerful, enabling the rapid convergence to local optima from random initial configurations.

Extending the optimization model from \citet{lai2024} to three dimensions introduces several significant challenges. While the dispersion penalty term is easy enough to compute in any finite number of dimensions, computing the boundary penalty term in three dimensions becomes computationally complex. First, \citet{lai2024} defines their boundary constraint term in a manner which heavily relies on planar geometry where counter-clockwise traversal of closed polygonal loops provides a canonical way to distinguish the interior from the exterior. We generalize boundary orientation of polygonal containers from two-dimensional vertex traversal to three dimensions by orienting the faces of the polyhedron using inward normals. By assigning an inward-pointing normal to each face, we explicitly define the interior half-space of the polyhedron. In our algorithm, these normals serve as the directional basis for calculating distance-to-boundary terms for dispersion points and performing feasibility checks.

Additionally, point-in-container testing is fundamental to the boundary constraint term, as it dictates whether to apply a repulsive force to maintain interiority or an attractive force to recover feasibility from the exterior. While point-in-polygon testing admits efficient deterministic routines \citep{wang2005}, no such efficient subroutines are devised for three-dimensional containers that scale well with the number of vertices. Consequently, we utilize a ray-casting-based point-in-polyhedron procedure \citep{horvat2012} together with the polyhedron's orientation to create projection-based distance calculations which handle complex or non-convex polyhedral containers. The ray-casting procedure draws multiple rays in a random directions from a given dispersion point and counts the number of faces it intersects. If a given ray intersects an even number of faces, the dispersion point is categorized as exterior; an odd number of intersections identifies the point as interior. Because rays may occasionally strike a vertex or edge leading to numerical ambiguity, we cast multiple rays and employ a majority-voting scheme to ensure a robust classification. These geometric tools provide the foundation needed to extend the underlying optimization model to three dimensions.

With this geometric framework in place, we now adapt the optimization model presented in \citet{lai2024} to three dimensions. The \textit{energy of a solution} $(c_1, \dots, c_p, D)$ will be measured by its quantification of constraint violation. Recall that high energy solutions have a high degree of constraint violation and low energy solutions have a low degree of constraint violation.

First, we define the penalty term for dispersion violation. Let 
\begin{equation}
	\ell_{i, j} = \max \set{0, 2D - d(c_i, c_j)},
	\label{eq:lij}
\end{equation}
hence we get $\ell_{i, j} = 0$ if and only if $d(c_i, c_j) \ge 2D$. Therefore the total pairwise dispersion penalty $\displaystyle{\sum_{i = 1}^{p - 1} \sum_{j = i + 1}^{p} \ell_{i, j}^2}$ equals $0$ if and only if all the dispersion points are sufficiently far from each other. By Rademacher's Theorem \citep{evans2015}, since $d$ is continuous, $\ell_{i, j}^2$ is differentiable almost-everywhere for all $i, j$.

Second, before we define the boundary constraint term, we need two auxiliary definitions which are three-dimensional analogues of the definitions presented by \citet{lai2024}:
    
    \textbf{Active Face:} Consider a polyhedral container $C$ with a face $F$ and a point $c$. Orient the normal vector $N_F$ of face $F$ inwards and denote $p \in F$ as the point in which $N_F$ originates from. \textit{A face $F$ is active with respect to $c$} if either $c \in C$ and $N_F \cdot (c - p) \ge 0$ or if $c \notin C$ and $N_F \cdot (c - p) < 0$.

\begin{figure}[h!]
    \centering
    \includegraphics[width=0.4\textwidth]{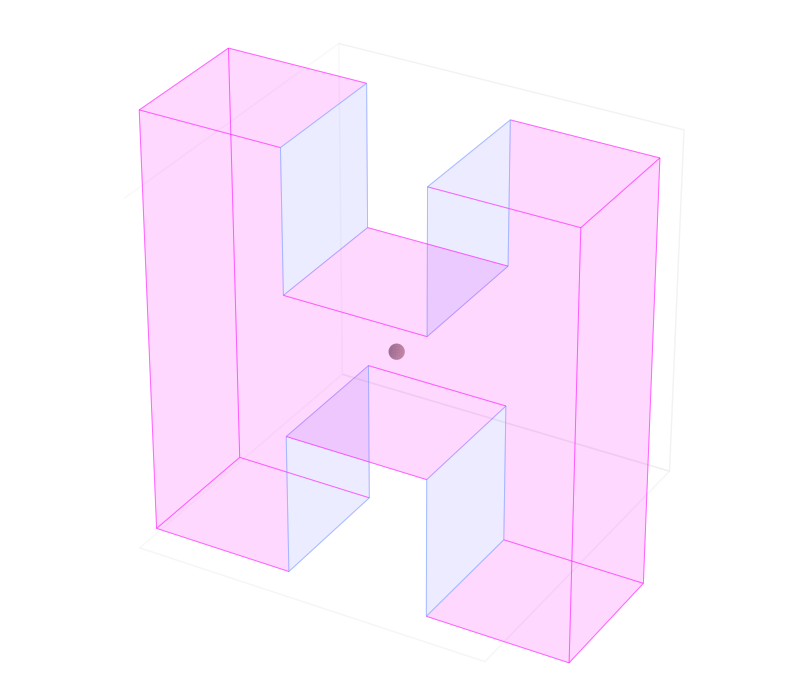}
    \caption{All active faces (magenta) of the H-box container for point $c = (1.5, 0.5, 1.5)$ inside the container.}
    \label{fig:active_planes_inside}
\end{figure}

\begin{figure}[h!]
    \centering
    \includegraphics[width=0.4\textwidth]{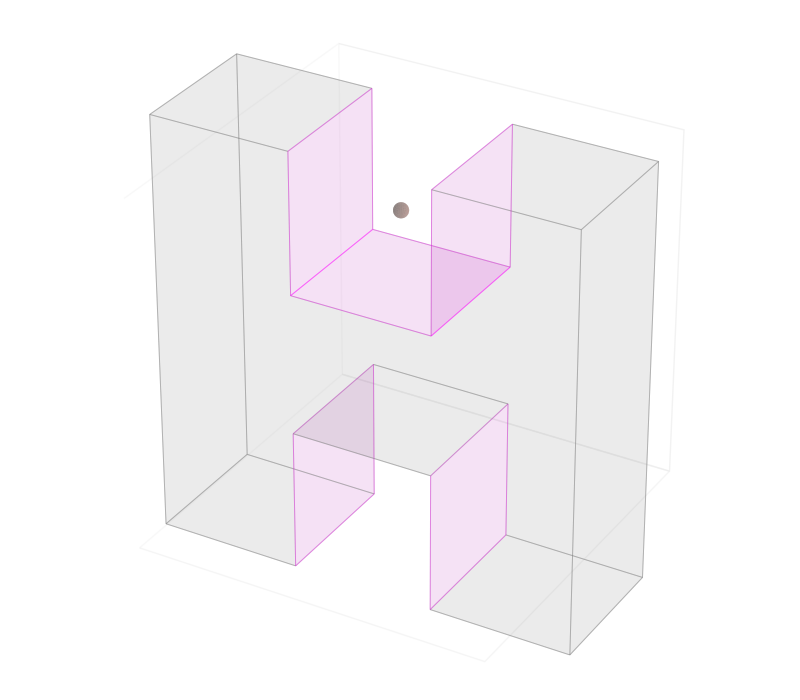}
    \caption{All active faces (magenta) of the H-box container for point $c = (1.5, 0.5, 2.5)$ outside the container.}
    \label{fig:active_planes_outside}
\end{figure}

Fig.~\ref{fig:active_planes_inside} depicts the active faces for a point inside the H-shaped container, and Fig.~\ref{fig:active_planes_outside} depicts the active faces for a point outside the H-shaped container. An active face can be thought of as an oriented plane that is geometrically relevant to the point $c$. If $c$ lies inside the polyhedron, then only the faces whose inward-normals minimally constrain how close $c$ is to the boundary of the container constitute the set of active faces. Alternatively, if $c$ lies outside the boundary of the polyhedron, an active face is blocking $c$ from entering the polyhedron. Active faces will contribute to the boundary constraint in the feasibility-residual function.

\textbf{Active Footpoint:} Consider a dispersion point $c$ and an active face $F$ of the container $C$. The \textit{foot} $h$ is the projection of $c$ onto the plane $P$ defined by the normal of $F$. We define $h$ to be an \textit{active footpoint} with respect to $F$ if $h$ lies on $F$. We denote $H^i$ to be the set of all active footpoints with respect to the dispersion point $c_i$.

\begin{figure}[h!]
    \centering
    \includegraphics[width=0.4\textwidth]{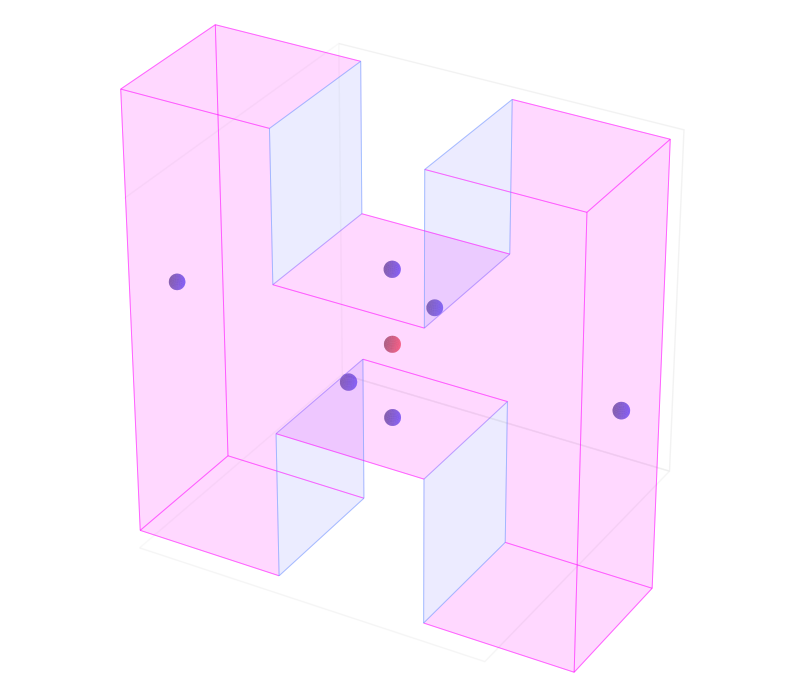}
    \caption{All active footpoints (blue) inside the H-box container with respect to (red) point $(1.5, 0.5, 1.5)$ inside container with activate faces highlighted in pink.}
    \label{fig:footpoints_inside}
\end{figure}

\begin{figure}[h!]
    \centering
    \includegraphics[width=0.4\textwidth]{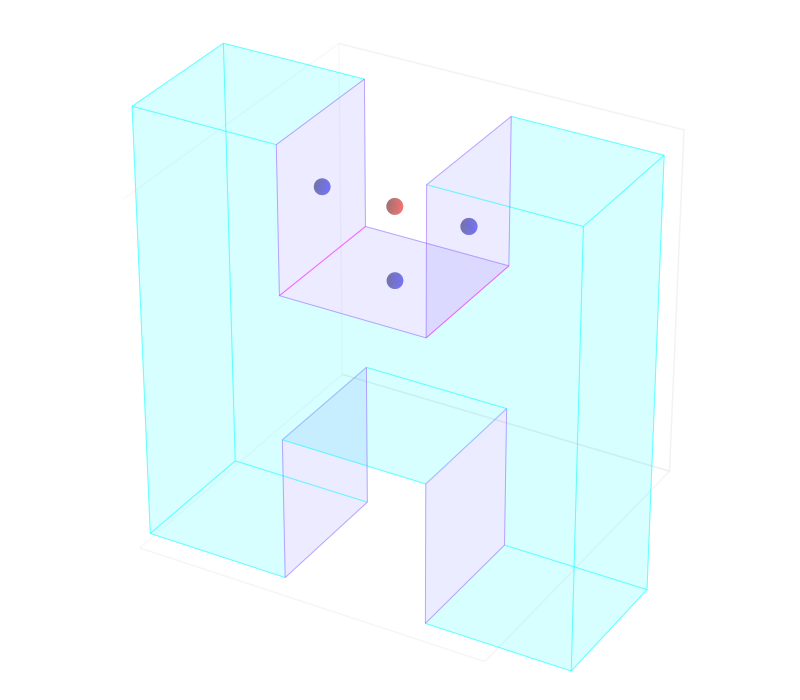}
    \caption{All active footpoints (blue) inside the H-box container with respect to (red) point $(1.5, 0.5, 2.5)$ outside container with active faces highlighted in purple.}
    \label{fig:footpoints_outside}
\end{figure}

Fig.~\ref{fig:footpoints_inside} depicts a case where the footpoints of a dispersion point are all active. Fig.~\ref{fig:footpoints_outside} depicts a case where not all footpoints of a dispersion point are active. The active footpoint $h$ of a dispersion point $c$ allows us to measure distances between dispersion points and the minimally relevant planes of the container, thereby providing a reference for the distance between a dispersion point and the container's boundary. When this perpendicular projection falls directly on the face $F$, $F$ is the boundary face that governs how close $c$ may move toward the container’s boundary, and the point-to-face distance is determined by the separation between $c$ and this footpoint.

Orientation equips the polyhedron with coherent inward normals, making the active face selection rule an orientation-consistent generalization of the counterclockwise active-edge mechanism used in the planar formulation \citep{lai2024}. Unlike the two-dimensional case, where the boundary is a closed, non-self-intersecting loop of segments with an inherent cyclic ordering, the boundary of a three-dimensional container forms a polyhedral complex with no analogous ordering. The extension from active edges \citep{lai2024} to active faces and their associated footpoints provides the geometric machinery needed to generalize boundary interactions from vertices and edges in two dimensions to vertices, edges, and faces in three dimensions. Thus, active faces constitute the minimal and geometrically natural extension of the active-edge concept from two to three dimensions, enabling a unified treatment of boundary constraints in arbitrary polyhedral domains.

Now, we define the boundary penalty for a dispersion point $c_i$ as
\begin{equation}
	O_i = \begin{cases}
	    D + 2\cdot (\min \set{d(c_i, v) : v \in V \cup H^i})^2 & c_i \notin C, \\ 
	    \displaystyle{\sum_{v \in V \cup E \cup H^i}} (\max\set{0, D - d(c_i, v)})^2 & c_i \in C.
	\end{cases}
	\label{eq:O_i}
\end{equation}
When $c_i \in C$, the boundary penalty is equal to the sum of all the squared distances between $c_i$  and the vertices, edges, and faces of the container in which the $c_i$ is at most $D$ away. In particular, given the open ball $B^o_D(c_i)$ of radius $D$ centered at $c_i$, each vertex, edge, and active footpoint of the container which has a non-empty intersection with $B^o_D(c_i)$ will contribute to the boundary penalty term $O_i$. If the distance to any vertex, edge, or footpoint from the dispersion point $c_i$ is larger than $D$, then no penalty is contributed. When $c_i \notin C$, a penalty is still imposed as the dispersion point is located outside the bounds of the container. In such a case, the minimum distance to a vertex or footpoint of the container is used as a basis for the penalty term. We add $D$ to this penalty term so that for points located outside the container, the boundary penalty term is always positive, and the penalty increases with $D$.

We now define the energy function $E: (\R^3)^p \times \R \to \R$ where $D > 0$ is fixed and $X = (c_1, \dots, c_p) \in (\R^3)^p$ by 
\begin{equation}
	E(X, D) = \sum_{i = 1}^{p - 1} \sum_{j = i + 1}^p \ell_{i, j}^2 + \sum_{i = 1}^p O_i.
	\label{eq:E}
\end{equation}

The first term of Equation~\ref{eq:E} deals with the dispersion constraint, or constraint (\text{ii}), of the CpDP. Recall, the first term of Equation~\ref{eq:E} is identically zero if and only if $d(c_i, c_j) \ge 2D$ for all distinct $i, j \in [p]$. When $D$ is fixed, we write $E_D(X) := E(X, D)$. Additionally, the second term deals with the containment constraint and the boundary constraint, or constraints (\textit{i}) and (\textit{iii}), respectively. Note, $O_i$ is zero if and only if $c_i \in C$ and $d(c_i, \partial C) \ge D$. Therefore, it is worth noting that $E_D(X)$ is identically zero if and only if for a fixed $D > 0$, the configuration $X = \set{c_1, \dots, c_p}$ is a feasible solution to the nonlinear optimization problem.

Equation~\ref{eq:E} gives a feasibility-residual function which lifts the optimization model proposed by \citet{lai2024} from two dimensions into three dimensions. Not only does it cleanly describe a feasibility-residual function for any polyhedral container, but is also differentiable almost-everywhere. The advantage of such a function allows one utilize gradient-based solvers (i.e., L-BFGS \citep{fletcher1987}) to find feasible solutions for a fixed $D \ge 0$ from initial random configurations.

\section{Global Optimization Algorithm}

\begin{figure}[h!] 
    \centering
    \includegraphics[scale=0.25]{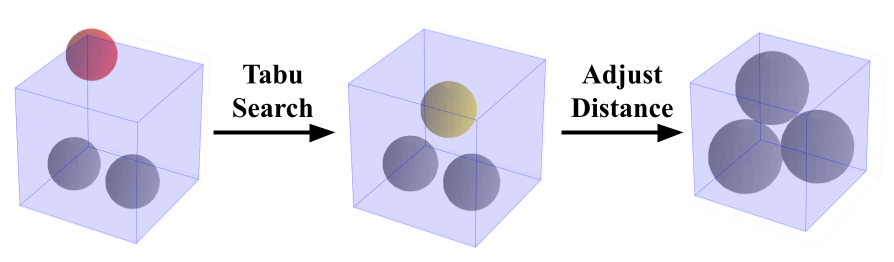}
    \caption{Overview of the global optimization algorithm. The process demonstrates a two-stage approach: first, Tabu Search is employed to relocate a candidate sphere (red) from an initial, infeasible position to a feasible position (yellow). Second, an Adjust Distance phase maximizes the uniform radius of all spheres to achieve an optimal packing configuration.}
    \label{fig:tabu_search_overview}
\end{figure}

In this section we propose our global optimization algorithm which builds upon the framework introduced by \citet{lai2024}. The global optimization process is structured as a multi-start pipeline to navigate the complex, non-convex solution space. In each global iteration, the algorithm restarts by generating an initial random configuration of $p$ points uniformly distributed within the polyhedral container. This multi-start approach is necessary because the problem’s non-smooth objective function and numerous distance constraints frequently trap standard optimization tools in poor local optima \citep{dai2021}. As shown in Fig.~\ref{fig:tabu_search_overview}, within each iteration, the Tabu Search procedure acts as a feasibility solver by minimizing the feasibility-residual function~\ref{eq:E} to find a valid configuration given a fixed minimum distance $D$. Once a feasible configuration is found, a distance adjustment procedure applies the Sequential Unconstrained Minimization Technique (SUMT) \citep{fiacco1964} to iteratively maximize $D$ while preserving the feasibility of the configuration.

\subsection{Tabu Search}

We employ a Tabu Search procedure \citep{glover1997}, following the framework used by \citet{lai2024}. Tabu Search is a combinatorial heuristic designed to escape poor local configurations (high-energy) and move the system toward a better (lower-energy) arrangement of dispersion points inside the container. The core principle behind the Tabu Search is the \textit{tabu list} which maintains a short term list of recent moves that are temporarily forbidden in subsequent iterations \citep{glover1997}. The tabu list ensures the algorithm does not fall into cycles of re-evaluating configurations. If the tabu list has reached its maximum allowable size, then the least recent addition to the tabu list is removed.

The Tabu Search heuristic incentivizes long-term exploration of the solution space. The maximum number of successive moves the algorithm is allowed to explore without improving from the initial configuration is specified by a parameter $\beta$. If Tabu Search makes $\beta$ iterations without improving on the last-improved configuration, then the procedure terminates. A small $\beta$ enables a quicker Tabu Search, but only scans for solutions which are a few moves different from the initial configurations. A large $\beta$ enables a deeper scan of neighboring solutions, but forces a longer runtime for Tabu Search. Tabu Search proves to be a powerful tool for efficiently scanning families of solutions using combinatorial optimization methods on a continuous region and transforming high-energy configurations to low-energy configurations.

For a given configuration $X = \set{c_1, \dots, c_p}$, we quantify the \textit{energy of a point} by the degree to which it violates the dispersion and boundary constraints. We define the energy of a dispersion point $c_i$ as
\begin{equation}
	E_D(c_i) = \sum_{j \ne i} \ell_{i, j}^2 + O_i
	\label{eq:dispersion_point_energy}
\end{equation}
where $\ell_{i, j}$ is the dispersion penalty term and $O_i$ captures the boundary violation term as in Equation~\ref{eq:O_i}. To propose new candidate positions for dispersion points, we consider \textit{vacancy sites}, which are random locations inside the container not currently occupied by a dispersion point. A vacancy site $c$ has energy
\begin{equation}
	E_D(c) = \sum_{i = 1}^p \ell_{c, j}^2 + O_{c}
	\label{eq:vacancy_energy}
\end{equation}
where $\ell_{c, j}$ replaces $c_i$ with $c$ in Equation~\ref{eq:lij} and $O_c$ replaces $c_i$ with $c$ in Equation~\ref{eq:O_i}. Vacancy sites with low energy represent positions that are more compatible with the current configuration, i.e., points where it would be the ``easiest'' to insert another ball.

A \textit{neighbor} $X_{nb}$ of a solution $X$ can be found by \textit{moving} a single dispersion point of $X$ to a new location. The \textit{neighborhood} $N(X)$ is the set of neighbors of $X$. To select the most promising neighbor of a given configuration $X$, we sample $N(X)$ by identifying the $Q$ highest-energy dispersion points and the $Q$ lowest-energy vacancy sites selected from a pool of $10p$ randomly generated candidates, where $p$ is the total number of dispersion points. While \citet{lai2024} suggest generating only $5p$ vacancy sites, we increased this sample size to $10p$ to maintain sufficient search breadth, compensating for the lower sampling density inherent in three-dimensional space compared to two. We analyze the possible moves between these sets to identify the pairing that offers the best potential reduction in energy.

For a specific dispersion-vacancy move, the resulting configuration is merely a rough approximation of a feasible solution. Since vacancy sites are generated randomly, simply placing a point at one of these coordinates rarely resolves the feasibility-residual value, defined in Equation~\ref{eq:E}, to zero. Therefore, we do not evaluate the quality of this raw move directly. Instead, we pass this intermediate configuration to our local optimization subroutine, \textsc{LocalOpt}. \textsc{LocalOpt} serves as the mandatory bridge between the discrete combinatorial search and the continuous solution space. It utilizes a gradient-based solver, specifically the Limited-memory Broyden-Fletcher-Goldfarb-Shanno (L-BFGS) algorithm \citep{fletcher1987}, to drive the rough candidate configuration toward the nearest smooth minimum of the energy function. Consequently, the Tabu Search algorithm never compares unrefined moves; it strictly compares the \textit{locally optimized} outcomes of those moves to decide if a step is an improvement.

However, even with this precise local refinement, the global solution space of Equation~\ref{eq:E} remains non-convex and populated with numerous local minima \citep{lai2024}. To prevent the algorithm from becoming trapped in the first deep valley it encounters, we integrate Monotonic Basin Hopping (\textsc{MBH}) as a heuristic for global exploration, as suggested by \citet{lai2024}. Unlike Tabu Search, which evolves a solution along a specific path, \textsc{MBH} samples the global landscape by ``hopping'' between independent basins (i.e., regions converging to a single minimum) \citep{back2013}. It operates by taking a solution $(X,D)$, performing random perturbations to jump into a neighboring basin, and then applying \textsc{LocalOpt}. The new local minimum is accepted only if it offers a strict improvement, ensuring the system can escape shallow local optima in search of a global minimum.

\begin{algorithm}[H]
\caption{Tabu Search}
\begin{algorithmic}[1]
    \Require $C$ (Container), $X$ (configuration), $D \in \R_{> 0}$, $\beta \in \Z_{> 0}$, $Q \in \Z_{> 0}$, $\epsilon \in \R_{> 0}$
    \State $Q \gets \min \set{Q, p}$
    \State $X \gets \textsc{LocalOpt}(X, D)$
    \State $X^{\text{best}} \gets X$, $\emph{NoImprove} \gets 0$, $\emph{TabuList} \gets \set{}$
    \While{$\emph{NoImprove} \le \beta$ and $E_D(X_{\text{best})} > \epsilon$} 
        \State $p_1, \dots, p_Q \gets \textsc{HighEnergyPoints}(X, D, Q)$ 
        \State $q_1, \dots, q_Q \gets \textsc{VacancySites}(X, D, Q)$
   
        \State $X^{\text{best}}_{\text{nb}} \gets \textsc{MovePoint}(X, p_1, q_1)$
        \For{$(i, j) \in [Q] \times [Q]$}
            \State $X_{\text{nb}} \gets \textsc{MovePoint}(X, p_i, q_j)$
            \State $X_{\text{nb}} \gets \textsc{LocalOpt}(X_{\text{nb}}, D)$
            \If {$\textsc{Move}(p_i, q_j)$ not forbidden}
                \State $X^{\text{best}}_{\text{nb}} \gets X_{\text{nb}}$ if $E_D(X_{\text{nb}}) < E_D(X^{\text{best}}_{\text{nb}})$
     
                \State $\emph{BestMove} \gets \textsc{Move}(p_i, q_j)$
            \EndIf
        \EndFor
        \State $X \gets X^{\text{best}}_{\text{nb}}$
        \State Update \emph{TabuList} with $\emph{BestMove}$
        \State $X \gets \textsc{MBH}(X, D)$
        \If{$E_D(X) < E_D(X^{\text{best}})$}
            \State $X^{\text{best}} \gets X$
     
            \State $\emph{NoImprove} \gets 0$
        \Else
            \State $\emph{NoImprove} \gets \emph{NoImprove} + 1$
        \EndIf
    \EndWhile
    \State \Return $X^{\text{best}}$
\end{algorithmic}
\label{alg:tabu_search}
\end{algorithm}

\begin{figure}[h!] 
    \centering
    \includegraphics[scale=0.2]{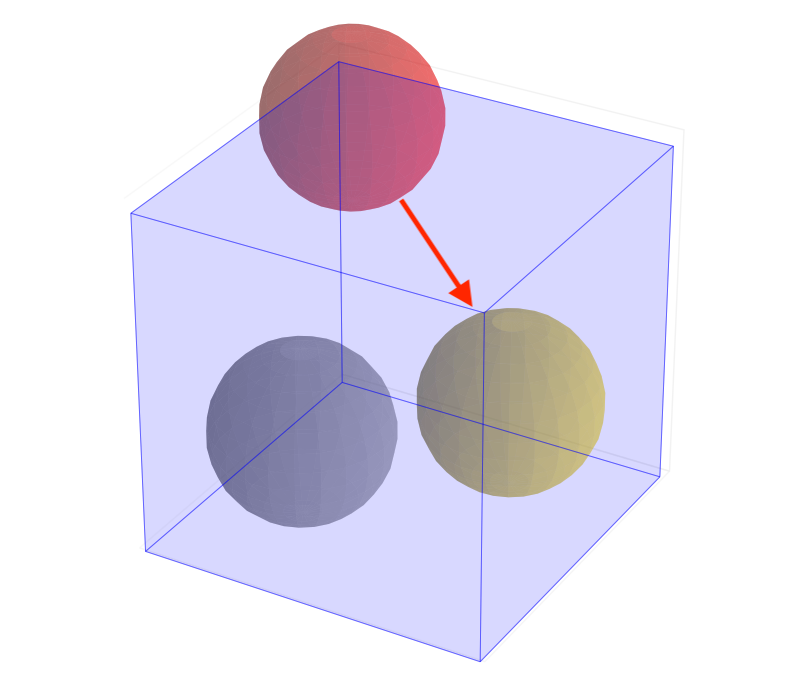}
    \caption{\textbf{Visualization of a candidate neighbor generation.} The red sphere represents the point selected for movement (e.g., a high-energy point). The red arrow indicates the proposed displacement vector, targeting the location currently occupied by the yellow sphere.}
    \label{fig:tabu_search}
\end{figure}

Algorithm~\ref{alg:tabu_search} outlines the Tabu Search procedure. To visualize the effect of Tabu Search, Fig.~\ref{fig:tabu_search} illustrates a candidate move for a fixed radius $D$ within a cubic container.

The dispersion point shown in red is in a high-energy state due to a boundary constraint violation (protrusion from the container). The yellow sphere represents a randomly generated vacancy site located in a feasible region. The red arrow indicates the proposed movement: displacing the violating particle to this vacancy site. In this specific example, because the target site satisfies all constraints (it is inside the boundary and non-overlapping), the feasibility-residual value is zero. Consequently, the subsequent \textsc{LocalOpt} routine will confirm this new configuration as stable with minimal. To prevent cycling, the algorithm explicitly forbids the reverse move in future iterations; the yellow dispersion point is temporarily restricted from returning to the high-energy neighborhood of the red sphere. This demonstrates the primary function of Tabu Search: executing coarse, global jumps to locate feasible regions, which are then refined by local optimization.

\subsection{Adjust Distance}

After the Tabu Search procedure is completed, the resulting configuration $X$ is a feasible solution, but typically not optimal. To find the nearest optimal solution $X^\star$ for the given configuration $X$, we employ the Sequential Unconstrained Minimization Technique (SUMT) \citep{fiacco1964} as proposed by \citet{lai2024} to convert the original CpDP (as a constrained optimization problem) into a sequence of unconstrained optimization problems and solve each until a feasible local minimum solution is found. Each unconstrained subproblem can be expressed as
\begin{equation}
	\textit{Minimize} \quad \Phi_\mu(X, D) = -D^2 + \mu \times E(X, D)
	\label{eq:phi_mu}
\end{equation} 
where $X$ is a candidate solution, $D$ is the radius of the solution $X$, $\mu$ is the penalty factor to define the unconstrained optimization problem.

Algorithm~\ref{alg:adjust_distance} is the method to adjust a feasible configuration from Algorithm~\ref{alg:tabu_search} into a local optimum. The distance adjustment method will iteratively minimize Equation~\ref{eq:phi_mu} for strictly increasing values of $\mu$. The $-D^2$ term in Equation~\ref{eq:phi_mu} is to ensure that when we minimize this function, we promote maximizing $D$. Equation~\ref{eq:O_i} contains $(\max\set{0, D - d(c_i, v)})^2$ which grows quadratically with $D$. The $-D^2$ term thus ensures Equation~\ref{eq:phi_mu} is not dominated by $E(X, D)$ \citep{lai2024}. The distance adjustment procedure slightly adjusts the configuration $X$ and the allowed minimum distance $D$ through local optimization to find a local optimum of the original CpDP.

\begin{algorithm}[H]
\caption{Adjust Distance}
\begin{algorithmic}[1]
    \Require $X_0$ (configuration), $D_0 \in \R_{> 0}$, $K \in \Z_{> 0}$
    \State $X \gets X_0$, $D \gets D_0$, $\mu \gets 10.0$
    \For{$i \gets 1$ to $K$}
        \State $(X, D) \gets \textsc{LocalOptPhi}(X, D, \mu)$
        \State $\mu \gets 5 \cdot \mu$
    \EndFor
    \State \Return $(X, D)$
\end{algorithmic}
\label{alg:adjust_distance}
\end{algorithm}

The \textsc{LocalOptPhi} function is a run of L-BFGS on the function $\Phi_\mu$ from an initial solution $(X_i, D_i)$ and penalty factor $\mu$. Empirical results suggest that $\mu$ can be initialized to at most $100$ so the gradient of $\Phi_\mu$ is not ill-conditioned at the initial iteration, allowing for a local optima to be located in the initial iteration. Note, we do not call \textsc{LocalOpt} in Algorithm~\ref{alg:adjust_distance}, since $D$ is fixed, hence minimizing $E_D(X)$ will not change $D$. Moreover, it would not help to call \textsc{LocalOptPhi} in Algorithm~\ref{alg:tabu_search} since Tabu Search is only designed for finding feasible configurations for a fixed $D$. It is important not to ``rush" finding a maximal value of $D$ as finding a stable and feasible configuration $X$ first often leads to a smoother maximization of $D$.

\subsection{Global Optimization}

The global optimization procedure is summarized in Algorithm~\ref{alg:global_opt}. This algorithm represents a three-dimensional analogue to the TSGO framework established by \citet{lai2024}. While \citet{lai2024} utilized a time-bounded loop with continuous dependency on the global best state, our approach adopts an iterative multi-start strategy. This modification decouples the search trajectories, allowing each iteration to explore a distinct portion of the solution space independently before comparing against the current global optimum.

By shifting to an iteration-based structure, we ensure that the algorithm extensively samples the configuration space. Each iteration begins with a fresh randomization, preventing the solver from becoming trapped in a previously found local optimum. This is particularly crucial for non-convex containers where the feasible region may be disconnected.

\begin{algorithm}[H]
\caption{Global Optimization}
\begin{algorithmic}[1]
	\Require $p \in \mathbb{N}$, $C$, $D_{\text{init}}$, $N$ (number of iterations)
	\Ensure Best configuration $(X^*, D^*)$
    
	\State $X \gets \textsc{RandomSolution}(p)$
	\State $X \gets \textsc{TabuSearch}(X, D)$
	\State $X, D \gets \textsc{AdjustDistance}(X, D, 15)$
	\State $X^*, D^* \gets X, D$
    
	\For {$i = 1$ \textbf{to} $N$}
		\State $D_i \gets D^*$
		\State $X_i \gets \textsc{RandomSolution}(p)$
		\State $X_i \gets \textsc{TabuSearch}(X_i, D_i)$
		\If{$E_{D_i}(X_i) < \varepsilon$}
			\State $(X_i, D_i) \gets \textsc{AdjustDistance}(X_i, D_i, 15)$
			\If{$D_i > D^*$}
	   			\State $X^*, D^* \gets X, D$
			\EndIf
		\EndIf
	\EndFor    
	\State \Return $(X^*, D^*)$
\end{algorithmic}
\label{alg:global_opt}
\end{algorithm}

Within each iteration, the algorithm retains the rigorous local search operators utilized by \citet{lai2024}. An initial configuration is generated via \textsc{RandomSolution}, followed by \textsc{TabuSearch} to locate the nearest local minimum. We then employ a strict feasibility check: only if the feasibility-residual value $E_D(X)$ for the configuration X falls below the tolerance threshold $\varepsilon$ does the algorithm proceed to \textsc{AdjustDistance}. This ensures that computational resources in the expensive refinement phase are reserved solely for valid, high-quality candidate solutions that have a genuine potential to improve upon the current global maximum $D^*$. At the local level, the optimization subroutines rely on a gradient-based solver governing Equation~\ref{eq:E}. This method exploits the specific geometry of the polyhedral container to achieve efficient convergence. Using gradient-based solvers over an almost-everywhere differentiable feasibility-residual function provides a smooth convergence path, offering greater stability than repulsion-based heuristics that rely on physical dynamical system simulations.

\section{Results and Evaluation}

Experimentation is carried out on a personal computer with an Apple M2 Pro 3.5 GHz processor and 32GB RAM running on a 64-bit ARM-based MacOS 26.2. CPU times are given in seconds. The algorithm is implemented in Python 3.11, and the source code is publicly available (\href{https://github.com/leadcatlab/continuous-p-dispersion}{Source Code}). Our specialized hybrid solver pairs Tabu Search with an L-BFGS optimizer to efficiently handle the problem's $3p+1$ variables and highly non-convex solution space, avoiding the local minima traps that stall standard generic solvers. We enforce a constraint violation tolerance of $\epsilon=10^{-8}$ and truncate reported results to eight decimal places, ensuring the reported configurations are strictly feasible within the container geometry.

We evaluate our algorithm on five polyhedral containers of increasing geometric complexity: the unit cube, unit-length tetrahedron, H-shaped box, star-shaped polyhedron, and a ribbed ventilation cage. The boundaries for each are defined using a vectorized representation of inward-pointing normals for efficient distance calculations. To evaluate robustness, we benchmarked our algorithm against Packomania optima \citep{packomania}, best-known results from \citet{kazakov2018}, and a computational baseline using the \texttt{scipy.optimize} \texttt{differential\_evolution} (DE) solver \citep{scipy2020}. The DE baseline was configured with the \texttt{best1bin} mutation strategy, a population size factor of 15, a 2400-second time limit, and a hard penalty for boundary constraints.

Algorithm parameters scaled with the number of dispersion points ($p$). Global iterations were set to 1 for small instances ($p < 10$), 3 for medium ($10 \le p < 30$), and 5 for large ($p \ge 30$). Tabu Search depth $\beta$ was similarly scaled to 5, 10, and 15, respectively, with a tabu-rearrangement parameter of $Q=3$ to align with \citet{lai2024}. The initial minimum distance $D_{\text{init}}$ was $0.01$ for small instances. For medium and large instances, $D_{\text{init}}$ was calculated using an initial packing density of $\rho_{\text{init}} = 0.3$ to ensure active penalty terms without over-congestion \citep{groemer1986, lai2024}:
\begin{equation}
    D_{\text{init}} = \sqrt[3]{\frac{3 \cdot V_{\text{container}} \cdot \rho_{\text{init}}}{4\pi p}}
    \label{eq:dinit_of_rho}
\end{equation}

\subsection{Cube}

We benchmark our algorithm in the unit cube against analytical optima for $p \le 10$ \citep{schaer1966} and best-known configurations for $p > 10$ from Packomania \citep{packomania}. We measure performance using absolute error ($A_{err} = D_{\mathrm{Ours}} - D_{\mathrm{Best}}$) and relative error ($R_{err} = A_{err}/D_{\mathrm{Best}} \times 100\%$).

\begin{figure}[h!]
    \centering
    \includegraphics[width=0.45\textwidth]{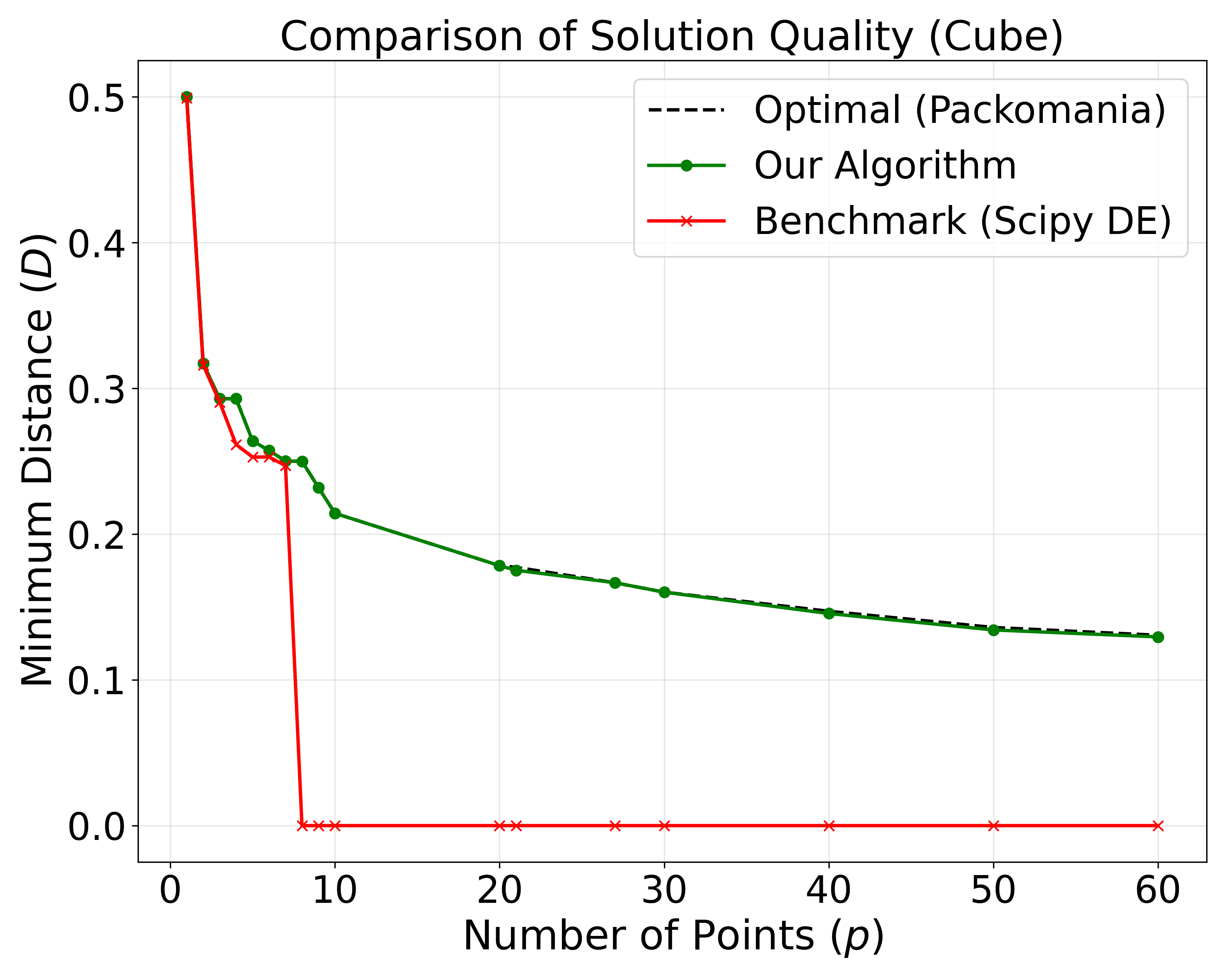}
    \caption{Performance comparison in the unit cube showing the best minimum distance $D$ for the proposed TSGO algorithm, the baseline DE solver, and the theoretical optima for $p$ points.}
    \label{fig:cube_comparison}
\end{figure}

Table~\ref{tab:cube_packomania_comparison} demonstrates high numerical accuracy. The algorithm consistently reaches the theoretical optimum for small instances ($p \le 10$, excluding $p=7$) and high-symmetry cases like $p=27$ ($<0.0001\%$ error). Despite the increasingly non-convex search space, $R_{err}$ for larger instances ($p=20$ to $60$) remains largely within $1\%$, with slight error increases at $p=21$ and $p=50$ likely indicating a high density of local optima. Furthermore, Fig.~\ref{fig:cube_comparison} confirms our method closely tracks theoretical values across all tested ranges, whereas the DE baseline's performance degrades rapidly as $p$ increases.

\begin{table}[h!]
    \centering
    \caption{Comparison of achieved minimum distances $D_{\mathrm{Ours}}$ with optimal values $D_{\mathrm{Best}}$ for $p=1, \dots, 60$ in the unit cube (Euclidean metric).}
    \label{tab:cube_packomania_comparison}
    \begin{tabular}{c|c|c|c|c}
        \hline
        $p$ & $D_{\mathrm{Ours}}$ & $D_{\mathrm{Best}}$ & $A_{err}$ & $R_{err} (\%)$ \\
        \hline
        1  & 0.50000000 & 0.50000000 & 0.00000000 & 0.0000\% \\
        2  & 0.31698729 & 0.31698730 & $-1.0\times10^{-8}$ & $-0.0000\%$ \\
        3  & 0.29289321 & 0.29289322 & $-1.0\times10^{-8}$ & $-0.0000\%$ \\
        4  & 0.29289321 & 0.29289322 & $-1.0\times10^{-8}$ & $-0.0000\%$ \\
        5  & 0.26393200 & 0.26393202 & $-2.0\times10^{-8}$ & $-0.0000\%$ \\
        6  & 0.25735926 & 0.25735931 & $-5.0\times10^{-8}$ & $-0.0000\%$ \\
        7  & 0.25009464 & 0.25013615 & $-4.2\times10^{-5}$ & $-0.0166\%$ \\
        8  & 0.24999999 & 0.25000000 & $-1.0\times10^{-8}$ & $-0.0000\%$ \\
        9  & 0.23205080 & 0.23205081 & $-1.0\times10^{-8}$ & $-0.0000\%$ \\
        10 & 0.21428565 & 0.21428571 & $-6.0\times10^{-8}$ & $-0.0000\%$ \\
        \hline
        20 & 0.17840693 & 0.17840720 & $-2.7\times10^{-7}$ & $-0.0002\%$ \\
        21 & 0.17521516 & 0.17721904 & $-2.0\times10^{-3}$ & $-1.1307\%$ \\
        27 & 0.16666660 & 0.16666667 & $-7.0\times10^{-8}$ & $-0.0000\%$ \\
        30 & 0.16018841 & 0.16018862 & $-2.1\times10^{-7}$ & $-0.0001\%$ \\
        40 & 0.14553682 & 0.14705882 & $-1.5\times10^{-3}$ & $-1.0350\%$ \\
        50 & 0.13426749 & 0.13595451 & $-1.7\times10^{-3}$ & $-1.2409\%$ \\
        60 & 0.12948983 & 0.13060194 & $-1.1\times10^{-3}$ & $-0.8515\%$ \\
        \hline
    \end{tabular}
\end{table}

\begin{figure}[h!]
    \centering
    \includegraphics[width=0.3\textwidth]{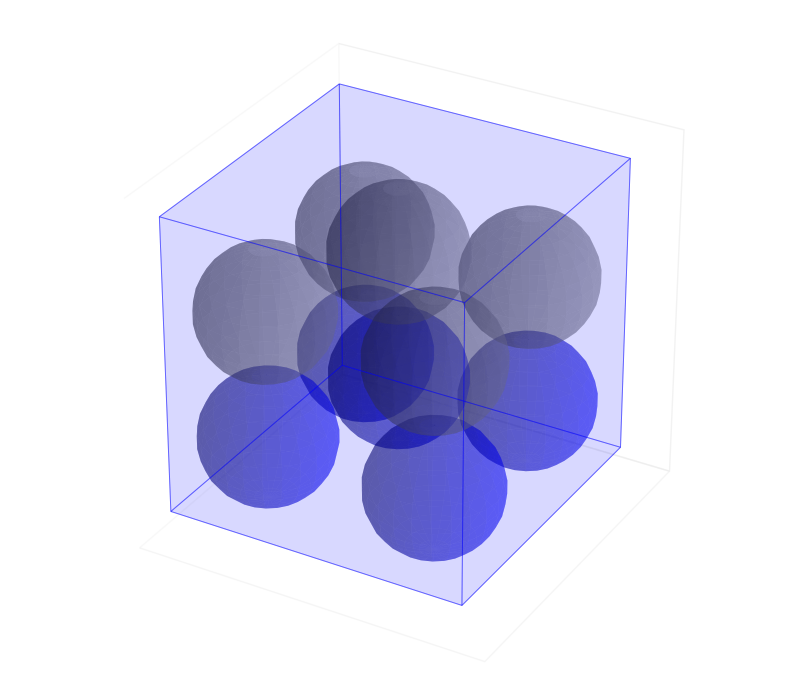}
    \caption{Best-achieved 10-point configuration for the $p$-dispersion problem in the unit cube using the Euclidean metric ($D = 0.214285$).}
    \label{fig:cube_10}
\end{figure}

\begin{figure}[h!]
    \centering
    \includegraphics[width=0.3\textwidth]{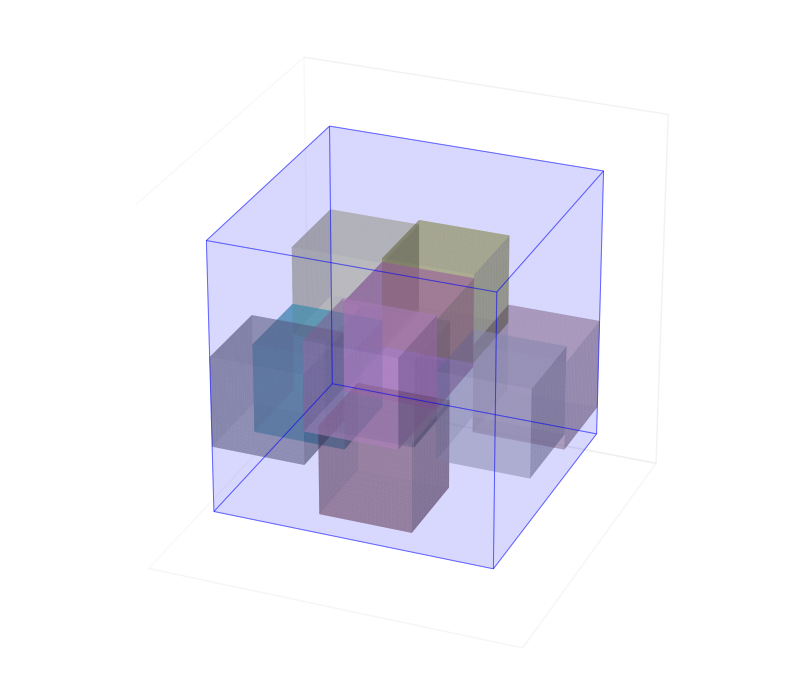}
    \caption{Best-achieved 10-point configuration for the $p$-dispersion problem in the unit cube using the $\ell_\infty$ metric ($D = 0.166666$).}
    \label{fig:cube_sup_10}
\end{figure}

To illustrate the algorithm's versatility across different geometries, we present visualizations in Fig.~\ref{fig:cube_10} and Fig.~\ref{fig:cube_sup_10} under the Euclidean and $\ell_\infty$ (sup-norm) metrics for packing ten points. The $\ell_\infty$ metric is particularly relevant as its geometry aligns with the rectilinear boundaries of the unit cube. As seen in Fig.~\ref{fig:cube_sup_10}, the algorithm effectively handles the flat-edged ``spheres'' of the sup-norm, achieving a configuration that is remarkably close to the theoretical optimum $(D^* = \frac16$).

\subsection{Tetrahedron}

We benchmark our algorithm against the best-known results for $p$-dispersion in the unit tetrahedron reported by \citet{kazakov2018} for $p = 1, \dots, 10$. Additionally, we report on our best-achieved minimum distances for $p = 20, 35, 56$ being the subsequent tetrahedral numbers. \citet{kazakov2018} mention that for three-dimensional space, the computational complexity of their algorithm increases significantly; hence they did not report on instances for $p > 10$. We execute our global optimization algorithm using the same parameters as in previous sections the packing density-based initialization ($\rho_{\text{init}} = 0.3$) for higher values of $p$. While \citet{kazakov2018} report results up to five decimal places, we report our results truncated to eight decimal places to maintain consistency with other sections of this paper.

Table~\ref{tab:tetra_kazakov} presents the absolute difference $A_{\text{err}} = D_{\mathrm{Ours}} - D_{\mathrm{Benchmark}}$ and the relative percentage improvement $R_{\text{err}}$ between our solution and the benchmark. As shown in the table, our algorithm consistently identifies packing configurations with equal or strictly larger minimum pairwise distances than \citet{kazakov2018}.

\begin{table}[h!]
    \centering
    \caption{Comparison of our achieved minimum distances ($D_{\mathrm{Ours}}$) against the best-known values ($D_{\mathrm{Benchmark}}$) from \citet{kazakov2018} for the unit tetrahedron.}
    \label{tab:tetra_kazakov}
    \begin{tabular}{c|c|c|c|c}
    \hline
    $p$ & $D_{\mathrm{Ours}}$ & $D_{\mathrm{Benchmark}}$ & $A_{\text{err}}$ & $R_{\text{err}}(\%)$ \\
    \hline
    1  & 0.20412415 & 0.20045 & 0.00367415 & +1.8330\% \\
    2  & 0.14494897 & 0.14204 & 0.00290897 & +2.0480\% \\
    3  & 0.14494896 & 0.14204 & 0.00290896 & +2.0480\% \\
    4  & 0.14494896 & 0.14204 & 0.00290896 & +2.0480\% \\
    5  & 0.12247447 & 0.12147 & 0.00100447 & +0.8269\% \\
    6  & 0.11387067 & 0.11243 & 0.00144067 & +1.2814\% \\
    7  & 0.11237234 & 0.11062 & 0.00175234 & +1.5841\% \\
    8  & 0.11237238 & 0.10600 & 0.00637238 & +6.0117\% \\
    9  & 0.11237238 & 0.10250 & 0.00987238 & +9.6316\% \\
    10 & 0.11237236 & 0.10200 & 0.01037236 & +10.1690\% \\
    \hline
    20 & 0.09175157 & -- & -- & -- \\
    35 & 0.07752526 & -- & -- & -- \\
    56 & 0.06190401 & -- & -- & -- \\
    \hline
    \end{tabular}
\end{table}

\begin{figure}[h!]
    \centering
    \includegraphics[width=0.45\textwidth]{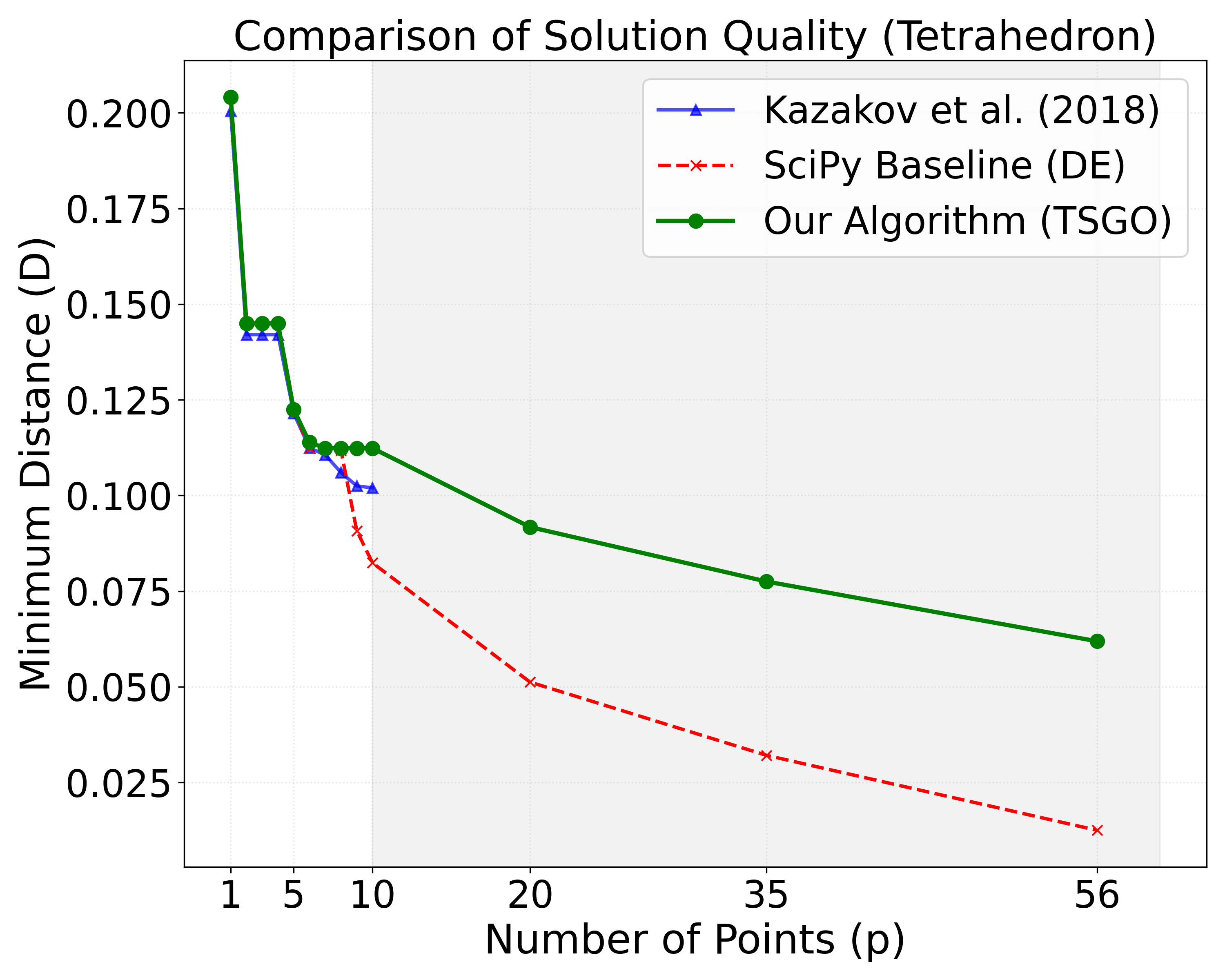}
    \caption{Performance comparison in the unit tetrahedron showing the best minimum distance $D$ for the proposed TSGO algorithm, the baseline DE solver, and the benchmarks from \citet{kazakov2018}.}
    \label{fig:tetra_gap}
\end{figure}

\begin{figure}[h!]
    \centering
    \begin{subfigure}[b]{0.45\textwidth}
        \centering
        \includegraphics[width=0.5\textwidth]{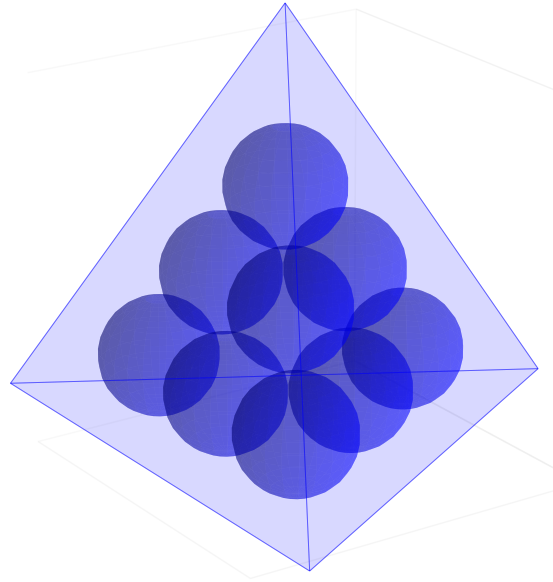}
        \caption{Our best-achieved configuration for $p=9$ ($D = 0.11237238$).}
        \label{fig:tetra_9_ours}
    \end{subfigure}
    \hfill
    \begin{subfigure}[b]{0.45\textwidth}
        \centering
        \includegraphics[width=0.6\textwidth]{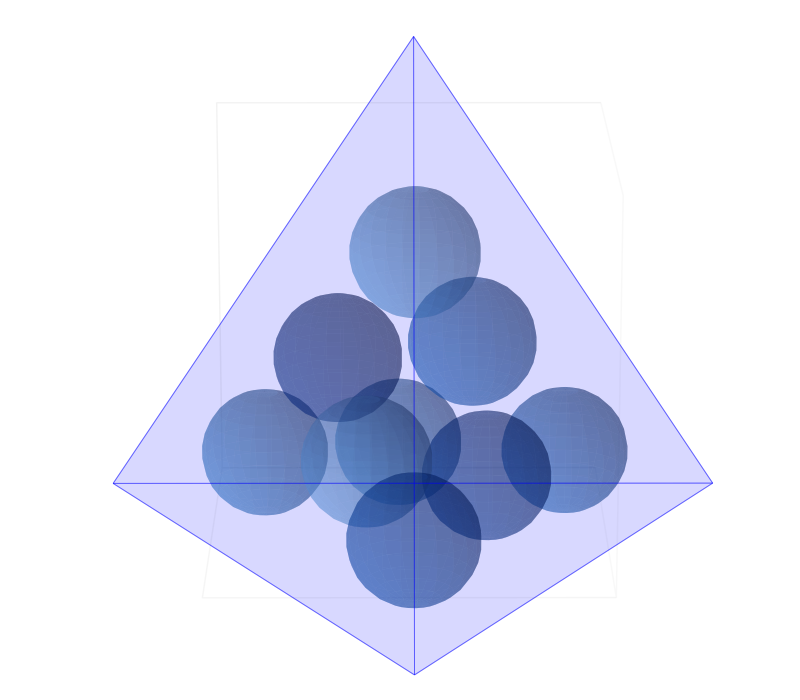}
        \caption{Feasible packing for $p=9$ using the radius ($D = 0.1025$) from \citet{kazakov2018}.}
        \label{fig:tetra_9_kazakov}
    \end{subfigure}
    \caption{Visual comparison of packing density for 9 dispersion points in the unit tetrahedron. Our method (a) achieves a significantly tighter packing than the benchmark parameters (b).}
    \label{fig:tetra_9}
\end{figure}

To further analyze these improvements, Fig.~\ref{fig:tetra_gap} illustrates the performance comparison between our algorithm, the results by \citet{kazakov2018}, and the DE solver. While the \texttt{scipy.optimize.differential\_evolution} baseline remains competitive for $p \le 8$, its performance degrades sharply for $p \ge 9$, with the error gap spiking significantly as the solver fails to navigate the increasingly complex configuration space. Fig.~\ref{fig:tetra_9} provides a visual comparison of the best configuration found by our algorithm against the corresponding configuration derived from the parameters in \citet{kazakov2018} for $p=9$.

The results demonstrate that our proposed method consistently outperforms both the existing literature and standard computational baselines. We achieve a minimum improvement of approximately 0.8\% for all tested values of $p$ compared to \citet{kazakov2018}, with the performance gap widening to over 10\% as $p$ reaches the limits of the previously reported benchmarks.

\subsection{H-Box}

Non-convex containers, such as the H-shaped polyhedron, present significant geometric challenges for the CpDP. Unlike convex domains, non-convex shapes often feature disconnected feasible regions and obstructed lines of sight, where the straight-line segment connecting two valid points may pass through the infeasible exterior \citep{dai2021}. These characteristics create numerous ``trapped'' local minima that can stall standard optimization algorithms \citep{lai2024}. The H-shaped polyhedron was created with inspiration from its two-dimensional analogue: the H-shaped polygon \citep{dai2021}.

To evaluate the performance of our algorithm, we compared our results against the SciPy DE solver. Table~\ref{tab:hbox_results} compares the best minimum distances ($D_{\mathrm{Ours}}$) achieved by our algorithm against this baseline ($D_{\mathrm{Baseline}}$) for $p = 1, \dots, 60$ points under the Euclidean metric. We report the absolute improvement $A_{\text{err}} = D_{\mathrm{Ours}} - D_{\mathrm{Baseline}}$ and the relative percentage improvement $R_{\text{err}}$.

\begin{table}[h!]
    \centering
    \caption{Comparison of achieved minimum distances for the H-shaped box ($D_{\mathrm{Ours}}$) against the Differential Evolution baseline ($D_{\mathrm{Baseline}}$).}
    \label{tab:hbox_results}
    \begin{tabular}{c|c|c|c|c}
    \hline
    $p$ & $D_{\mathrm{Ours}}$ & $D_{\mathrm{Baseline}}$ & $A_{\text{err}}$ & $R_{\text{err}}(\%)$ \\
    \hline
    1  & 0.50000001 & 0.49999742 & +0.00000259 & +0.001\% \\
    2  & 0.50000000 & 0.49999183 & +0.00000817 & +0.002\% \\
    3  & 0.50000000 & 0.49997360 & +0.00002640 & +0.005\% \\
    4  & 0.50000000 & 0.49997731 & +0.00002269 & +0.005\% \\
    5  & 0.50000000 & 0.49996936 & +0.00003064 & +0.006\% \\
    6  & 0.50000000 & 0.49938900 & +0.00061100 & +0.122\% \\
    7  & 0.50000000 & 0.47443668 & +0.02556332 & +5.388\% \\
    8  & 0.41093592 & 0.27026219 & +0.14067373 & +52.051\% \\
    9  & 0.39725847 & 0.36077084 & +0.03648763 & +10.114\% \\
    10 & 0.39360889 & 0.22369906 & +0.16990983 & +75.955\% \\
    \hline
    20 & 0.30427141 & 0.00000000 & +0.30427141 & $\infty$ \\
    30 & 0.28046438 & 0.00000000 & +0.28046438 & $\infty$ \\
    40 & 0.26696214 & 0.00000000 & +0.26696214 & $\infty$ \\
    50 & 0.24122813 & 0.00000000 & +0.24122813 & $\infty$ \\
    56 & 0.24999983 & 0.00000000 & +0.24999983 & $\infty$ \\
    60 & 0.19986383 & 0.00000000 & +0.19986383 & $\infty$ \\
    \hline
    \end{tabular}
\end{table}

The results demonstrate a stark performance gap between the two methods. While the baseline solver remains competitive for $p \le 7$, it fails to report a nontrivial configuration for larger values of $p$, configurations in which the minimum distance is $D = 0.0$. In contrast, our proposed algorithm successfully identifies high-quality solutions across the entire range, maintaining structural feasibility even at high packing densities.

\begin{figure}[h!]
    \centering
    \includegraphics[width=0.4\textwidth]{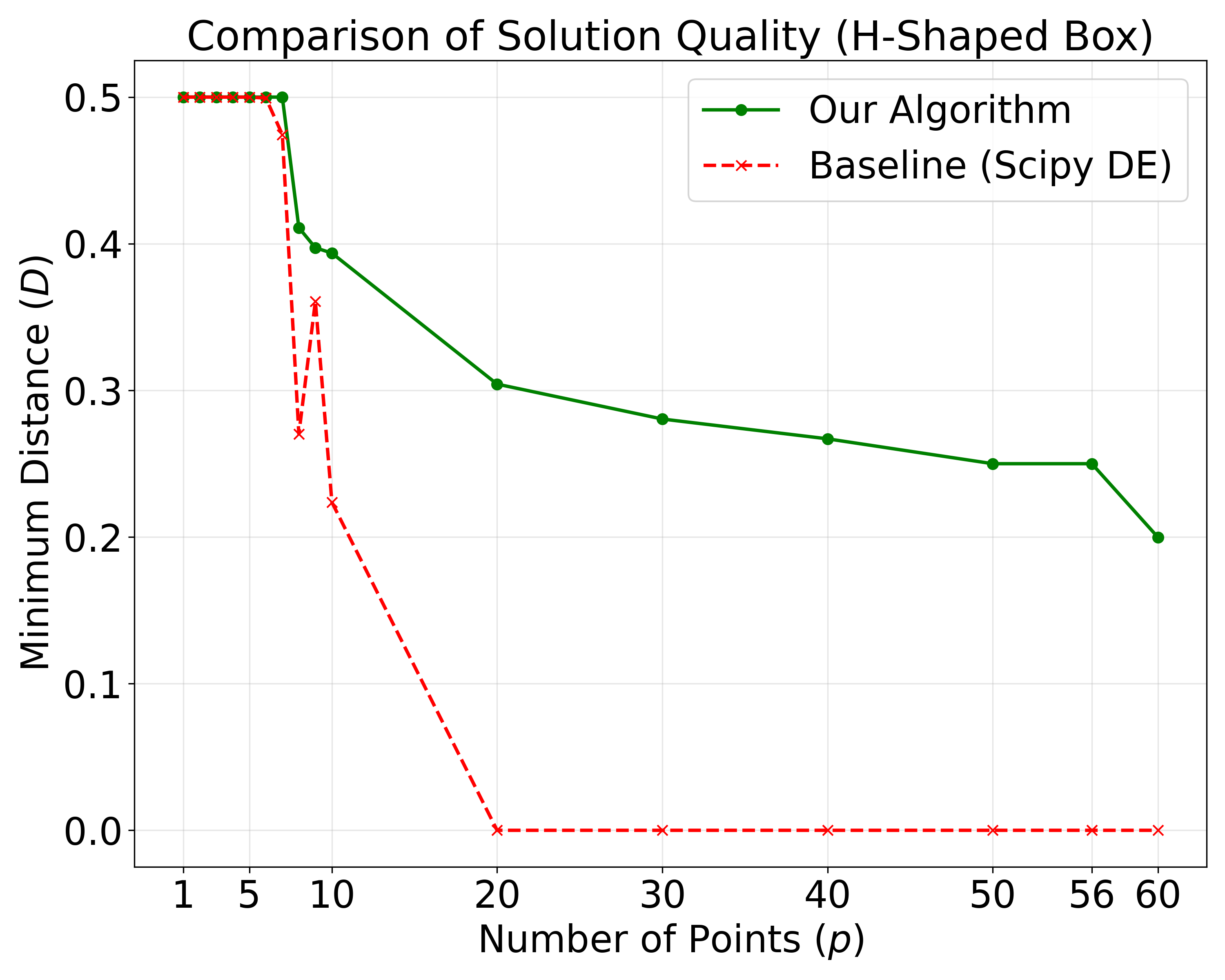}
    \caption{Performance comparison on the H-shaped box. The baseline solver (red) fails completely for $p \ge 20$ ($D=0$), whereas our algorithm (green) effectively explores the non-convex passage to maintain non-trivial separation distances.}
    \label{fig:hbox_comparison}
\end{figure}

To validate the performance of our method in these non-convex domains, Fig.~\ref{fig:hbox_comparison} illustrates the best minimum distance ($D$) achieved by each solver. The divergence highlights the baseline's inability to navigate the narrow passages and internal corners of the H-box, while our algorithm remains robust.

\begin{figure}[h!]
    \centering
    \includegraphics[width=0.35\textwidth]{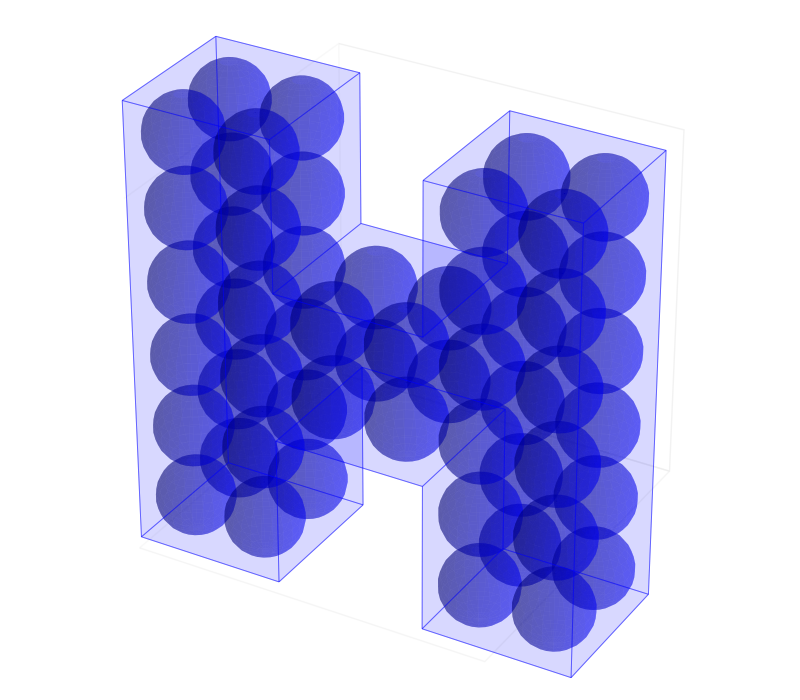}
    \caption{Best achieved configuration of $p=56$ dispersion points within the non-convex H-shaped container ($D = 0.24999983$). The visualization highlights the algorithm's ability to resolve tight packing constraints within narrow, non-convex passages.}
    \label{fig:h_box_56}
\end{figure}

Fig.~\ref{fig:h_box_56} visualizes the best-achieved configuration for $p=56$ points. This high-density packing illustrates the algorithm's ability to maintain structural stability across non-convex junctions. Notably, the spheres form a coordinated ``bridge'' through the central passage while maintaining a rigid structure in the lateral segments, effectively utilizing the restricted volume to achieve $D \approx 0.25$. This configuration confirms that our Tabu Search method effectively navigates the disconnected feasible regions that typically obstruct global optimizers.

\subsection{Star-Shaped Polyhedron}

The star-shaped polyhedron draws inspiration from the Kepler-Poinsot polyhedron \citep{coxeter1989} and serves as a rigorous test case for non-convex optimization due to its distinctive central convex core connected to six narrow, tapering spikes. This originally developed geometry creates significant ``bottlenecks'' between the feasible sub-regions, making it difficult for optimization algorithms to migrate points from one spike to another. Unlike convex containers where local descent directions are often globally informative, the star-shaped polyhedron requires an algorithm capable of traversing infeasible space to escape the local optima found within individual spikes.

In Table~\ref{tab:pointy_star_results}, we report the best minimum distances ($D_{\mathrm{Ours}}$) achieved by our algorithm for a range of $p$ values. To demonstrate the scalability of our approach on highly constrained non-convex domains, we report results for large-scale instances up to $p=100$.

\begin{table}[h!]
    \centering
    \caption{Achieved minimum distances for the star-shaped polyhedron ($D_{\mathrm{Ours}}$).}
    \label{tab:pointy_star_results}
    \begin{tabular}{c|c||c|c}
    \hline
    $p$ & $D_{\mathrm{Ours}}$ & $p$ & $D_{\mathrm{Ours}}$ \\
    \hline
    1   & 1.41421356 & 9   & 0.72337216 \\
    2   & 0.97597051 & 10  & 0.69809740 \\
    3   & 0.90296401 & 20  & 0.57376175 \\
    4   & 0.90296402 & 30  & 0.49246630 \\
    5   & 0.90296400 & 40  & 0.45625711 \\
    6   & 0.90296396 & 50  & 0.39072828 \\
    7   & 0.81657905 & 60  & 0.36665003 \\
    8   & 0.76614274 & 100 & 0.33677304 \\
    \hline
    \end{tabular}
\end{table}

The results indicate that our algorithm effectively exploits the available volume even as the number of points increases. For $p=1$, the algorithm identifies the globally optimal placement within the central core. As $p$ increases to 6, the algorithm successfully distributes points into the six available spikes, maintaining a high separation distance of $D \approx 0.9029$. Beyond $p=7$, crowding effects force a gradual reduction in the minimum distance $D$ as the central core and spikes become filled.

\begin{figure}[h!]
    \centering
    \includegraphics[width=0.5\textwidth]{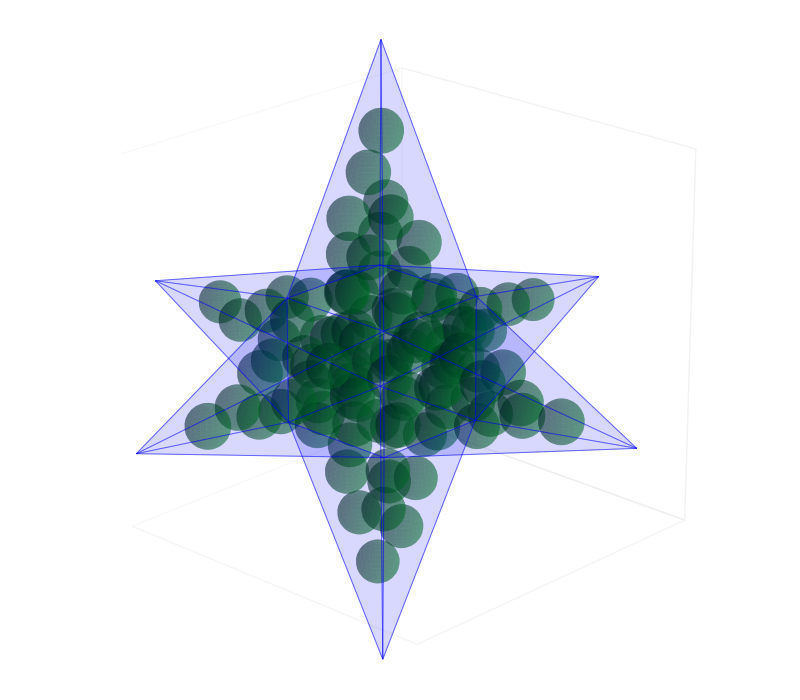}
    \caption{Best achieved configuration of $p=100$ dispersion points within the star-shaped container ($D = 0.33677304$).}
    \label{fig:ps_100}
\end{figure}

Fig.~\ref{fig:ps_100} visualizes the best-achieved configuration for $p=100$. This result is particularly notable as it demonstrates the algorithm's capability to pack dispersion points uniformly into the tapering extremities of the spikes. Such regions are typically inaccessible to standard gradient-based solvers due to the narrow feasible widths and the high likelihood of points becoming trapped in local optima near the spike tips. The ability of our Tabu Search hybrid to relocate points across the non-convex "bottlenecks" is evidenced by the high packing density achieved in all six spikes simultaneously.

\subsection{Ribbed Ventilation Cage}

Finally, we apply our algorithm to the ``ribbed ventilation cage,'' a highly complex non-convex polyhedron designed to model modern industrial computer casings and heat sinks. Unlike the H-shaped box or Star-shaped polyhedron, which feature large distinct sub-regions, this container is characterized by numerous corridors and separated legs. The inclusion of numerous internal pillars and ribs creates a ``maze-like'' feasible region, introducing frequent discontinuities that pose severe challenges for traditional gradient-based nonlinear optimizers.

The particular design of this ribbed ventilation cage polyhedral container is original, although the complexity is practically motivated. In thermal engineering, components utilize ribbed structures, fins, and microstructures to maximize surface area and improve heat dissipation capabilities \citep{xin2025}. Solving the CpDP on such a geometry is therefore relevant for optimizing the layout of sensors, cooling pipes, or modular components within densely packed electronic enclosures. We evaluate our algorithm using both the Euclidean metric ($\ell_2$), representing standard spherical component packing, and the sup-norm ($\ell_\infty$), which is critical for modeling the placement of cubic or box-shaped components often found in hardware design.

\begin{figure}[H]
	\centering
    \begin{minipage}[t]{0.4\textwidth}
        \centering
        \includegraphics[width=\textwidth]{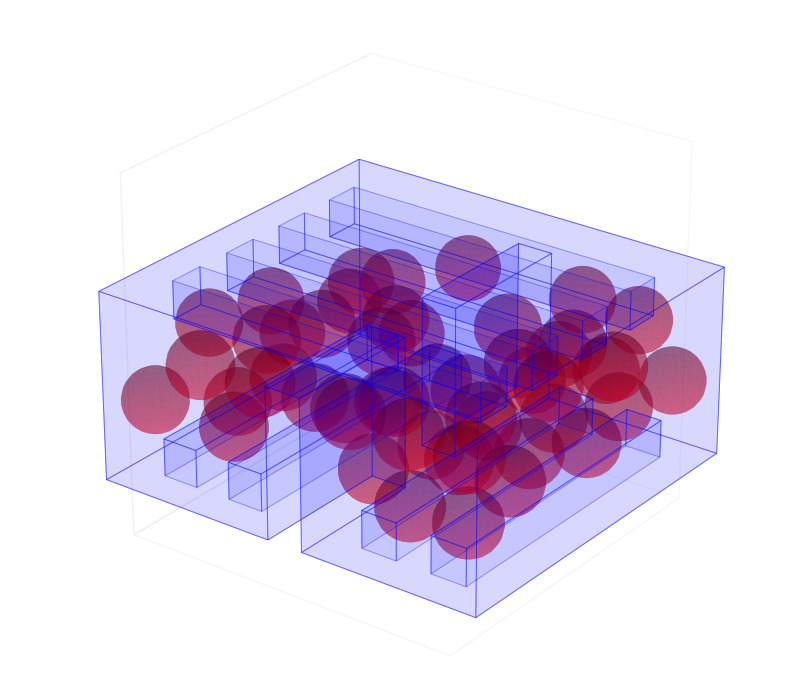}
        \caption{Best achieved configuration for $p = 50$ points using the Euclidean metric ($D = 0.8247$).}
        \label{fig:cage_50_l2}
    \end{minipage}\hfill
    \begin{minipage}[t]{0.4\textwidth}
        \centering
        \includegraphics[width=\textwidth]{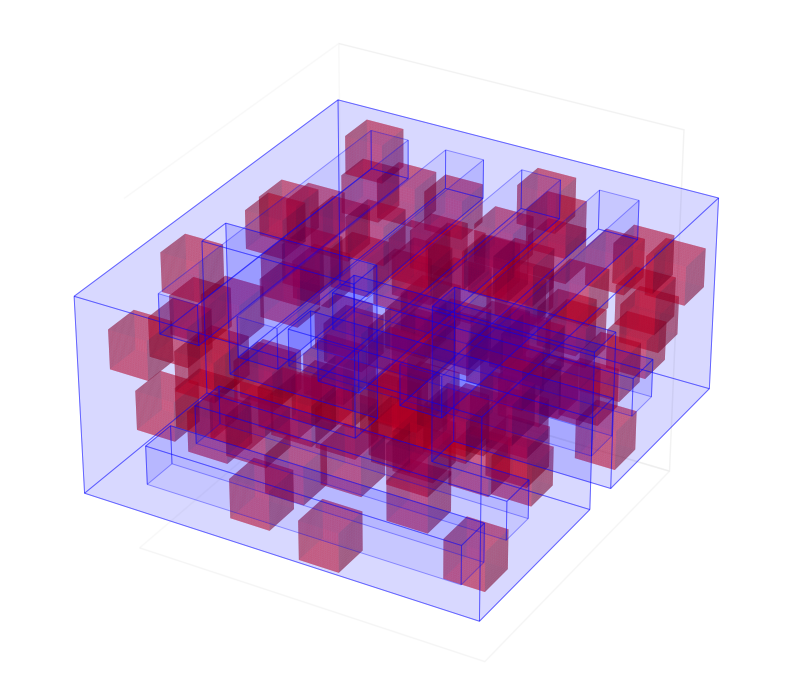}
        \caption{Best achieved configuration for $p=100$ points using the $\ell_\infty$ metric ($D = 0.4999$).}
        \label{fig:cage_100_sup}
    \end{minipage}
\end{figure}

Fig.~\ref{fig:cage_50_l2} and Fig.~\ref{fig:cage_100_sup} visualize the results of our algorithm for dispersing $p = 50$ points with respect to the Euclidean matric and $p = 100$ points with respect to the $\ell_\infty$ metric within this highly constrained environment.

The successful dispersion of $p=50$ points in this environment demonstrates the robustness of our approach. Despite the numerous internal obstacles that fracture the search space, the algorithm achieves a uniform distribution, effectively utilizing the narrow corridors between pillars. This result validates the algorithm's capability to handle large-scale dispersion problems ($p \ge 50$) within containers that exhibit both macro-scale non-convexity and micro-scale geometric obstructions.

\subsection{Computational Efficiency and Scalability}

\begin{figure}[H]
    \centering
    \includegraphics[width=0.4\textwidth]{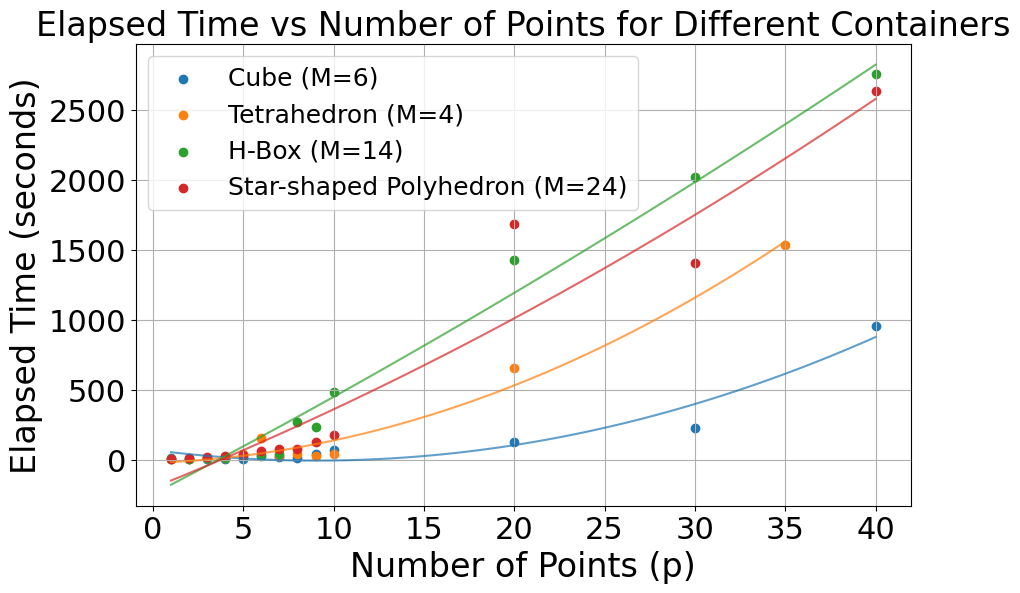}
    \caption{\textbf{Comparative runtime analysis.} Runtime (seconds) versus dispersion points $p$. Solid lines indicate quadratic regression fits.}
    \label{fig:combined_runtime}
\end{figure}

To evaluate scalability, we recorded the wall-clock runtime required to converge to a feasible solution ($\varepsilon = 10^{-8}$). Fig.~\ref{fig:combined_runtime} compares runtimes across four container geometries (excluding the ribbed cage, where complex boundary checks overwhelmingly dominate). 

The empirical data closely matches quadratic regression fits, aligning with the $\mathcal{O}(p^2)$ complexity of pairwise distance checks and $\mathcal{O}(p \cdot M)$ complexity of boundary checks, where $M$ is the number of faces. For small instances ($p \le 10$), convergence takes under 50 seconds—vastly outperforming the 2400-second Differential Evolution baseline, which frequently timed out on non-convex shapes.

As $p$ increases, geometric complexity dictates the scaling coefficient. Simple shapes like the unit tetrahedron ($M=4$) and cube ($M=6$) exhibit shallow growth. Conversely, the star-shaped polyhedron ($M=24$) scales more steeply due to severe bottlenecks and increased ``Active Face'' evaluations. Crucially, all geometries maintain polynomial scaling, avoiding the exponential explosion common in combinatorial approaches. Even for the highly complex star-shaped container at $p=100$, the algorithm converges in a tractable $\sim$4,900 seconds.

\section{Conclusion}

This paper introduces a robust global optimization framework for the Continuous $p$-Dispersion Problem (CpDP) in three-dimensional polyhedral containers under any metric. By defining a differentiable feasibility-residual function and extending the Tabu Search algorithm \cite{lai2024}, we successfully transition from restrictive two-dimensional planar checks to a comprehensive three-dimensional formulation utilizing a ray-casting procedure and projection-based distance calculations to faces, edges, and vertices. Empirical evaluations confirm that our method matches analytical optima in standard shapes \cite{schaer1966} and establishes new high-dispersion benchmarks for complex, non-convex polyhedra \cite{kazakov2018}, offering a powerful tool for real-world applications such as robotic collision avoidance and facility layout \cite{abdulghafoor2021, loganathan2024}.

Because our energy formulation relies on dimension-independent vector operations, it provides a natural foundation for extension into $n$-dimensional spaces ($\mathbb{R}^n$). The primary adaptation required for higher dimensions is generalizing the geometric container representation, specifically determining how to precisely model $(n-1)$-dimensional bounding hyperplanes and extending the ray-casting validity checks. Furthermore, future work will focus on incorporating more sophisticated, non-planar geometries such as curved surfaces. While polyhedral distances currently utilize efficient orthogonal projections, curved boundaries require representations like spline patches or implicit surfaces, necessitating the development of computationally efficient distance evaluation methods that do not rely on closed-form projections.

\pagebreak

\begin{widetext}
\section*{Acknowledgements}
The authors would like to thank Ram Padmanabhan for their valuable feedback and editorial assistance throughout the preparation of this manuscript. 
\end{widetext}

\bibliographystyle{unsrtnat} 
\bibliography{continuous-p-dispersion-in-three-dimensions}

\end{document}